\documentclass[a4paper,12pt]{article}
\usepackage[utf8]{inputenc}
\usepackage[T1]{fontenc}
\usepackage{amsfonts,amsmath,amssymb,amsthm}
\theoremstyle{definition}

\newtheorem{conjecture}{Conjecture}
\theoremstyle{plain}

\newtheorem{theorem}{Theorem}
\usepackage{graphicx}
\usepackage{caption}
\usepackage{pgf,tikz,circuitikz,subfig}
\usetikzlibrary{arrows,shapes,positioning}
\usepackage{booktabs,multirow,tabularx}
\usepackage{hyperref}
\hypersetup{colorlinks,%
	citecolor={red}, 
	urlcolor={blue},
	linkcolor={blue},
	breaklinks={true},
	pagebackref={true},
	hyperindex={true},
}
\usepackage{url}
\usepackage{cite}
\usepackage[left=2cm,right=2cm,top=2cm,bottom=2cm]{geometry}
\usepackage[ruled]{algorithm}
\usepackage{algorithmic}
\usepackage{newtxtext}
\usepackage{courier}
\usepackage{enumitem}
\newlist{abbrv}{itemize}{1}
\setlist[abbrv,1]{label=,labelwidth=0.9in,align=parleft,noitemsep,leftmargin=!}

\newcommand{\ub}[1]{\overline{#1}}

\newcommand{\geo}[1]{\mathtt{#1}}

\title{Tight bounds on the maximal area of small polygons: Improved Mossinghoff polygons}
\author{Christian Bingane\thanks{D\'{e}partement de math\'{e}matiques de g\'{e}nie industriel, Polytechnique Montr\'{e}al, Montreal, QC, Canada. Email: \url{christian.bingane@polymtl.ca}}}
\begin{document}
\maketitle
\begin{abstract}
A small polygon is a polygon of unit diameter. The maximal area of a small polygon with $n=2m$ vertices is not known when $m \ge 7$. In this paper, we construct, for each $n=2m$ and $m\ge 3$, a small $n$-gon whose area is the maximal value of a one-variable function. We show that, for all even $n\ge 6$, the area obtained improves by $O(1/n^5)$ that of the best prior small $n$-gon constructed by Mossinghoff. In particular, for $n=6$, the small $6$-gon constructed has maximal area.
\end{abstract}

\section{Introduction}
The {\em diameter} of a polygon is the largest Euclidean distance between pairs of its vertices. A polygon is said to be {\em small} if its diameter equals one. For an integer $n \ge 3$, the maximal area problem consists in finding a small $n$-gon with the largest area. The problem was first investigated by Reinhardt~\cite{reinhardt1922} in 1922. He proved that
\begin{itemize}
	\item for all $n \ge 3$, the value $\frac{n}{2}\sin \frac{\pi}{n} - \frac{n}{2}\tan \frac{\pi}{2n}$ is an upper bound on the area of a small $n$-gon;
	\item when $n$ is odd, the regular small $n$-gon is the unique optimal solution;
	\item when $n=4$, there are infinitely many optimal solutions, including the small square;
	\item when $n \ge 6$ is even, the regular small $n$-gon is not optimal.
\end{itemize}

\begin{figure}
	\centering
	\subfloat[$(\geo{R}_6,0.649519)$]{
		\begin{tikzpicture}[scale=4]
			\draw[dashed] (0,0) -- (0.4330,0.2500) -- (0.4330,0.7500) -- (0,1) -- (-0.4330,0.7500) -- (-0.4330,0.2500) -- cycle;
			\draw (0,0) -- (0,1);
			\draw (0.4330,0.2500) -- (-0.4330,0.7500);
			\draw (0.4330,0.7500) -- (-0.4330,0.2500);
		\end{tikzpicture}
	}
	\subfloat[$(\geo{R}_5^+,0.672288)$]{
		\begin{tikzpicture}[scale=4]
			\draw[dashed] (0,0) -- (0.5000,0.3633) -- (0.3090,0.9511) -- (0,1) -- (-0.3090,0.9511) -- (-0.5000,0.3633) -- cycle;
			\draw (0,1) -- (0,0) -- (0.3090,0.9511) -- (-0.5000,0.3633) -- (0.5000,0.3633) -- (-0.3090,0.9511) -- (0,0);
		\end{tikzpicture}
	}
	\subfloat[$(\geo{P}_6^*,0.674981)$]{
		\begin{tikzpicture}[scale=4]
			\draw[dashed] (0,0) -- (0.5000,0.4024) -- (0.3438,0.9391) -- (0,1) -- (-0.3438,0.9391) -- (-0.5000,0.4024) -- cycle;
			\draw (0,1) -- (0,0) -- (0.3438,0.9391) -- (-0.5000,0.4024) -- (0.5000,0.4024) -- (-0.3438,0.9391) -- (0,0);
		\end{tikzpicture}
		\label{figure:6gon:U6}
	}
	\caption{Three small $6$-gons $(\geo{P}_6,A(\geo{P}_6))$}
	\label{figure:6gon}
\end{figure}
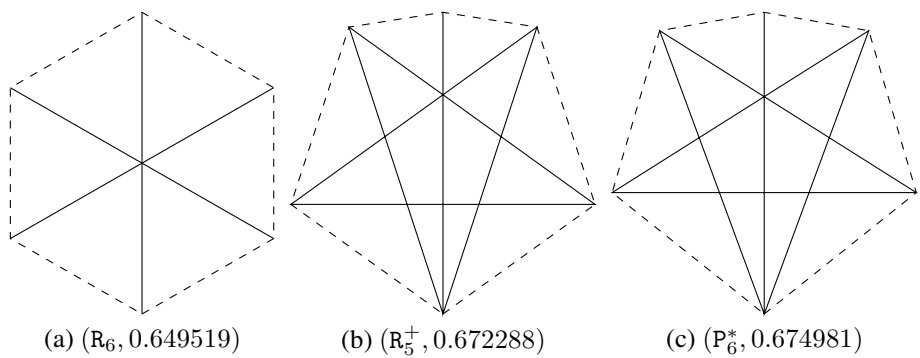

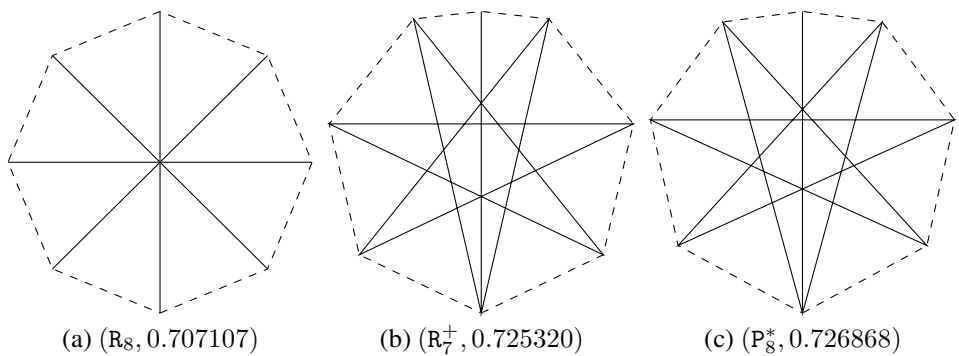
\begin{figure}
	\centering
	\subfloat[$(\geo{R}_8,0.707107)$]{
		\begin{tikzpicture}[scale=4]
			\draw[dashed] (0,0) -- (0.3536,0.1464) -- (0.5000,0.5000) -- (0.3536,0.8536) -- (0,1) -- (-0.3536,0.8536) -- (-0.5000,0.5000) -- (-0.3536,0.1464) -- cycle;
			\draw (0,0) -- (0,1);
			\draw (0.3536,0.1464) -- (-0.3536,0.8536);
			\draw (0.5000,0.5000) -- (-0.5000,0.5000);
			\draw (0.3536,0.8536) -- (-0.3536,0.1464);
		\end{tikzpicture}
	}
	\subfloat[$(\geo{R}_7^+,0.725320)$]{
		\begin{tikzpicture}[scale=4]
			\draw[dashed] (0,0) -- (0.4010,0.1931) -- (0.5000,0.6270) -- (0.2225,0.9749) -- (0,1) -- (-0.2225,0.9749) -- (-0.5000,0.6270) -- (-0.4010,0.1931) -- cycle;
			\draw (0,1) -- (0,0) -- (0.2225,0.9749) -- (-0.4010,0.1931) -- (0.5000,0.6270) -- (-0.5000,0.6270) -- (0.4010,0.1931) -- (-0.2225,0.9749) -- (0,0);
		\end{tikzpicture}
	}
	\subfloat[$(\geo{P}_8^*,0.726868)$]{
		\begin{tikzpicture}[scale=4]
			\draw[dashed] (0,0) -- (0.4091,0.2238) -- (0.5000,0.6404) -- (0.2621,0.9650) -- (0,1) -- (-0.2621,0.9650) -- (-0.5000,0.6404) -- (-0.4091,0.2238) -- cycle;
			\draw (0,1) -- (0,0) -- (0.2621,0.9650) -- (-0.4091,0.2238) -- (0.5000,0.6404) -- (-0.5000,0.6404) -- (0.4091,0.2238) -- (-0.2621,0.9650) -- (0,0);
		\end{tikzpicture}
		\label{figure:8gon:U8}
	}
	\caption{Three small $8$-gons $(\geo{P}_8,A(\geo{P}_8))$}
	\label{figure:8gon}
\end{figure}

When $n \ge 6$ is even, the maximal area problem is solved for $n \le 12$. The case $n=6$ was solved by Bieri~\cite{bieri1961} in 1961 and Graham~\cite{graham1975} in 1975, the case $n=8$ by Audet, Hansen, Messine, and Xiong~\cite{audet2002} in 2002, and the cases $n=10$ and $n=12$ by Henrion and Messine~\cite{henrion2013} in 2013. Both optimal $6$-gon and $8$-gon are represented in Figure~\ref{figure:6gon:U6} and Figure~\ref{figure:8gon:U8}, respectively. In 2017, Audet~\cite{audet2017} showed that the regular small polygon has the maximal area among all equilateral small polygons.

The diameter graph of a small polygon is the graph with the vertices of the polygon, and an edge between two vertices exists only if the distance between these vertices equals one. Diameter graphs of some small polygons are shown in Figure~\ref{figure:4gon}, Figure~\ref{figure:6gon}, and Figure~\ref{figure:8gon}. The solid lines illustrate pairs of vertices which are unit distance apart. In 2007, Foster and Szabo~\cite{foster2007} proved that, for even $n \ge 6$, the diameter graph of a small $n$-gon with maximal area has a cycle of length $n-1$ and one additional edge from the remaining vertex. From this result, they provided a tighter upper bound on the maximal area  of a small $n$-gon when $n \ge 6$ is even.

\begin{figure}
	\centering
	\subfloat[$(\geo{R}_4,0.5)$]{
		\begin{tikzpicture}[scale=4]
			\draw[dashed] (0,0) -- (0.5000,0.5000) -- (0,1) -- (-0.5000,0.5000) -- cycle;
			\draw (0,0) -- (0,1);
			\draw (0.5000,0.5000) -- (-0.5000,0.5000);
		\end{tikzpicture}
	}
	\subfloat[$(\geo{R}_3^+,0.5)$]{
		\begin{tikzpicture}[scale=4]
			\draw[dashed] (0.5000,0.8660) -- (0,1) -- (-0.5000,0.8660);
			\draw (0,1) -- (0,0) -- (0.5000,0.8660) -- (-0.5000,0.8660) -- (0,0);
		\end{tikzpicture}
	}
	\caption{Two small $4$-gons $(\geo{P}_4,A(\geo{P}_4))$}
	\label{figure:4gon}
\end{figure}
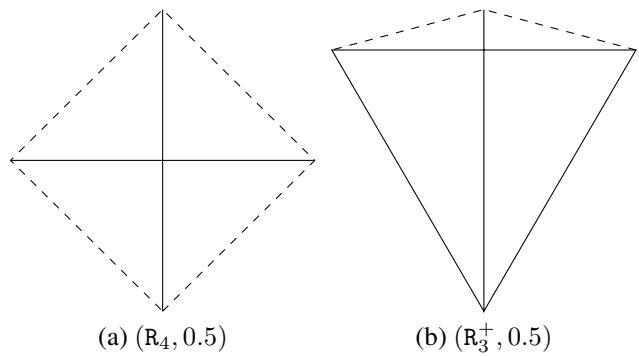

For even $n\ge 8$, exact solutions in the maximal area problem appear to be presently out of reach. However, tight lower bounds on the maximal area can be obtained analytically. For instance, Mossinghoff~\cite{mossinghoff2006b} constructed a family of small $n$-gons, for even $n\ge 6$, and proved that the areas obtained cannot be improved for large $n$ by more than $c/n^3$, for a certain positive constant $c$. By contrast, the areas of the regular small $n$-gons cannot be improved for large $n$ by more than $\pi^3/(16n^2)$ when $n \ge 6$ is even. In this paper, we propose tighter lower bounds on the maximal area of small $n$-gons when $n \ge 6$ is even. Thus, the main result of this paper is the following:

\begin{theorem}\label{thm:Bn}
	Suppose $n = 2m$ with integer $m \ge 3$. Let $\ub{A}_n := \frac{n}{2}\sin \frac{\pi}{n} - \frac{n-1}{2}\tan \frac{\pi}{2n-2}$ denote an upper bound on the area $A(\geo{P}_n)$ of a small $n$-gon $\geo{P}_n$~\cite{foster2007}. Let $\geo{M}_n$ denote the small $n$-gon constructed by Mossinghoff~\cite{mossinghoff2006b} for the maximal area problem. Then there exists a small $n$-gon $\geo{B}_n$ such that
	\[
	\ub{A}_n - A(\geo{B}_n) = \frac{(5303-456\sqrt{114})\pi^3}{5808n^3} + O\left(\frac{1}{n^4}\right) < \frac{3\pi^3}{40n^3} + O\left(\frac{1}{n^4}\right)
	\]
	and
	\[
	A(\geo{B}_n) - A(\geo{M}_n) = \frac{3d\pi^3}{n^5} + O\left(\frac{1}{n^6}\right)
	\]
	with
	\[
	\begin{aligned}
		d &= \frac{25\pi^2(1747646-22523\sqrt{114})}{4691093528} + \frac{32717202988-3004706459\sqrt{114}}{29464719680}\\
		&+ (-1)^{\frac{n}{2}} \frac{15\pi(10124777-919131\sqrt{114})}{852926096}\\
		&=
		\begin{cases}
			0.0836582354\ldots &\text{if $n \equiv 2 \bmod 4$,}\\
			0.1180393778\ldots &\text{if $n \equiv 0 \bmod 4$.}
		\end{cases}
	\end{aligned}
	\]
	Moreover, $\geo{B}_6$ is the largest small $6$-gon.
\end{theorem}

The remainder of this paper is organized as follows. Section~\ref{sec:ngon} recalls principal results on the maximal area problem. We prove Theorem~\ref{thm:Bn} in Section~\ref{sec:Bn}. We conclude the paper in Section~\ref{sec:conclusion}.

\section{Areas of small polygons}\label{sec:ngon}
Let $A(\geo{P})$ denote the area of a polygon $\geo{P}$. Let $\geo{R}_n$ denote the regular small $n$-gon. We have
\[
A(\geo{R}_n) =
\begin{cases}
	\frac{n}{2}\sin \frac{\pi}{n} - \frac{n}{2}\tan \frac{\pi}{2n} &\text{if $n$ is odd,}\\
	\frac{n}{8}\sin \frac{2\pi}{n} &\text{if $n$ is even.}\\
\end{cases}
\]
For all even $n\ge 6$, $A(\geo{R}_n) < A(\geo{R}_{n-1})$~\cite{audet2009}. This suggests that $\geo{R}_n$ does not have maximum area for any even $n\ge 6$. Indeed, when $n$ is even, we can construct a small $n$-gon with a larger area than $\geo{R}_n$ by adding a vertex at distance $1$ along the mediatrix of an angle in $\geo{R}_{n-1}$. We denote this $n$-gon by $\geo{R}_{n-1}^+$ and we have
\[
A(\geo{R}_{n-1}^+) = A(\geo{R}_{n-1}) + \sin \frac{\pi}{2n-2} - \frac{1}{2}\sin \frac{\pi}{n-1}.
\]

\begin{theorem}[Reinhardt~\cite{reinhardt1922}, Foster and Szabo~\cite{foster2007}]\label{thm:area:opt}
	For all $n \ge 3$, let $A_n^*$ denote the maximal area among all small $n$-gons.
	\begin{itemize}
		\item When $n$ is odd, $A_n^* = \frac{n}{2}\sin \frac{\pi}{n} - \frac{n}{2}\tan \frac{\pi}{2n}$ is only achieved by $\geo{R}_n$.
		\item $A_4^* = 1/2$ is achieved by infinitely many $4$-gons, including $\geo{R}_4$ and~$\geo{R}_3^+$ illustrated in Figure~\ref{figure:4gon}.
		\item When $n\ge 6$ is even, the diameter graph of an optimal $n$-gon has a cycle of length $n-1$ plus one additional edge from the remaining vertex and $A_n^* < \ub{A}_n := \frac{n}{2} \sin \frac{\pi}{n} - \frac{n-1}{2} \tan \frac{\pi}{2n-2}$.
	\end{itemize}
\end{theorem}

When $n\ge 6$ is even, the maximal area~$A_n^*$ is known for $n \le 12$. Bieri~\cite{bieri1961} and Graham~\cite{graham1975} determined analytically that $A_6^* = 0.674981\ldots > A(\geo{R}_{5}^+)$, and this value  is only achieved by the small $6$-gon shown in Figure~\ref{figure:6gon:U6}. Audet, Hansen, Messine, and Xiong~\cite{audet2002} proved that $A_8^* = 0.726868\ldots > A(\geo{R}_{7}^+)$, which is only achieved by the small $8$-gon represented in Figure~\ref{figure:8gon:U8}. Henrion and Messine~\cite{henrion2013} found that $A_{10}^* = 0.749137\ldots > A(\geo{R}_{9}^+)$ and $A_{12}^* = 0.760729\ldots > A(\geo{R}_{11}^+)$.

\begin{conjecture}
	\label{thm:area:sym}
	For even $n \ge 6$, an optimal $n$-gon has an axis of symmetry corresponding to the pendant edge in its diameter graph.
\end{conjecture}

From Theorem~\ref{thm:area:opt}, we note that $\geo{R}_{n-1}^+$ has the optimal diameter graph. Conjecture~\ref{thm:area:sym} is only proven for $n=6$ and this is due to Yuan~\cite{yuan2004}. However, the largest small polygons obtained by~\cite{audet2002} and~\cite{henrion2013} are a further evidence that the conjecture may be true.

\begin{figure}
	\centering
	\subfloat[$(\geo{M}_6,0.673186)$]{
		\begin{tikzpicture}[scale=4]
			\draw[dashed] (0,0) -- (0.5000,0.4362) -- (0.3701,0.9290) -- (0,1) -- (-0.3701,0.9290) -- (-0.5000,0.4362) -- cycle;
			\draw (0,1) -- (0,0) -- (0.3701,0.9290) -- (-0.5000,0.4362) -- (0.5000,0.4362) -- (-0.3701,0.9290) -- (0,0);
		\end{tikzpicture}
	}
	\subfloat[$(\geo{M}_{8},0.725976)$]{
		\begin{tikzpicture}[scale=4]
			\draw[dashed] (0,0) -- (0.3988,0.2265) -- (0.5000,0.6649) -- (0.2813,0.9596) -- (0,1) -- (-0.2813,0.9596) -- (-0.5000,0.6649) -- (-0.3988,0.2265) -- cycle;
			\draw (0,1) -- (0,0) -- (0.2813,0.9596) -- (-0.3988,0.2265) -- (0.5000,0.6649) -- (-0.5000,0.6649) -- (0.3988,0.2265) -- (-0.2813,0.9596) -- (0,0);
		\end{tikzpicture}
	}
	\subfloat[$(\geo{M}_{10},0.749029)$]{
		\begin{tikzpicture}[scale=4]
			\draw[dashed] (0,0) -- (0.3310,0.1396) -- (0.5000,0.4454) -- (0.4463,0.7687) -- (0.2167,0.9762) -- (0,1) -- (-0.2167,0.9762)  -- (-0.4463,0.7687) -- (-0.5000,0.4454) -- (-0.3310,0.1396) -- cycle;
			\draw (0,1) -- (0,0) -- (0.2167,0.9762) -- (-0.3310,0.1396) -- (0.4463,0.7687) -- (-0.5000,0.4454) -- (0.5000,0.4454) -- (-0.4463,0.7687) -- (0.3310,0.1396) -- (-0.2167,0.9762) -- (0,0);
		\end{tikzpicture}
	}
	\caption{Mossinghoff polygons $(\geo{M}_n,A(\geo{M}_n))$}
	\label{figure:Mn}
\end{figure}
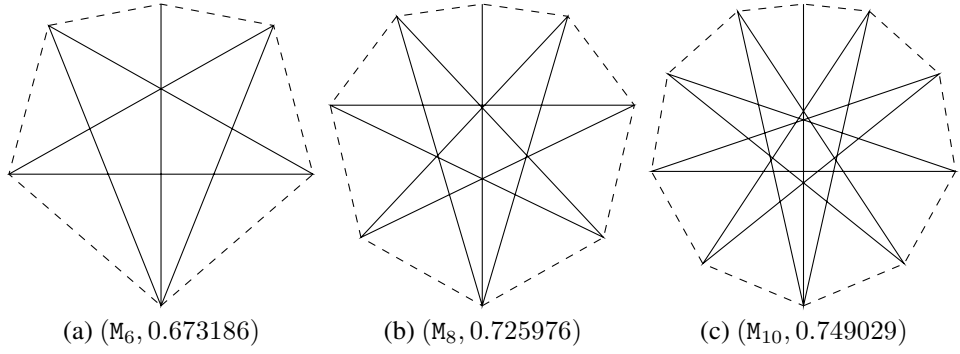

For even $n\ge 6$, $\geo{R}_{n-1}^+$ does not provide the tightest lower bound for $A_n^*$. Indeed, Mossinghoff~\cite{mossinghoff2006b} constructed a family of small $n$-gons $\geo{M}_n$, illustrated in Figure~\ref{figure:Mn}, such that
\[
\ub{A}_n - A(\geo{M}_n) = \frac{(5303-456\sqrt{114})\pi^3}{5808n^3} + O\left(\frac{1}{n^4}\right) < \frac{3\pi^3}{40n^3} + O\left(\frac{1}{n^4}\right)
\]
for all even $n \ge 6$. On the other hand,
\[
\begin{aligned}
	\ub{A}_n - A(\geo{R}_n) &= \frac{\pi^3}{16n^2} + O\left(\frac{1}{n^3}\right),\\
	\ub{A}_n - A(\geo{R}_{n-1}^+) &= \frac{5\pi^3}{48n^3} + O\left(\frac{1}{n^4}\right)
\end{aligned}
\]
for all even $n\ge 6$. In the next section, we propose a tighter lower bound for $A_n^*$.

\section{Proof of Theorem~\ref{thm:Bn}}\label{sec:Bn}
\begin{figure}
	\centering
	\begin{tikzpicture}[scale=8]
		\draw[dashed] (0,0) node[below]{$\geo{v}_0(0,0)$} -- (0.4068,0.2215) node[right]{$\geo{v}_5(x_5,y_5)$} -- (0.5000,0.6432) node[right]{$\geo{v}_3(x_3,y_3)$} -- (0.2619,0.9651) node[right]{$\geo{v}_1(x_1,y_1)$} -- (0,1) node[above]{$\geo{v}_7(0,1)$} -- (-0.2619,0.9651) node[left]{$\geo{v}_6(x_6,y_6)$} -- (-0.5000,0.6432) node[left]{$\geo{v}_4(x_4,y_4)$} -- (-0.4068,0.2215) node[left]{$\geo{v}_2(x_2,y_2)$} -- cycle;
		\draw (0,1) -- (0,0) -- (0.2619,0.9651) -- (-0.4068,0.2215) -- (0.5000,0.6432) -- (-0.5000,0.6432) -- (0.4068,0.2215) -- (-0.2619,0.9651) -- (0,0);
		\draw (0.0655,0.2413) arc (74.82:90.00:0.25) node[midway,above]{$\alpha_0$};
		\draw (0.0946,0.7793) arc (228.01:254.82:0.25) node[midway,below]{$\alpha_1$};
		\draw (-0.1801,0.3269) arc (24.93:48.01:0.25) node[midway,right]{$\alpha_2$};
		\draw (0.2500,0.6432) arc (180.00:204.93:0.25) node[midway,left]{$\alpha_3$};
	\end{tikzpicture}
	\caption{Definition of variables $\alpha_0, \alpha_1, \ldots, \alpha_{\frac{n}{2}-1}$: Case of $n=8$ vertices}
	\label{figure:model:optimal}
\end{figure}
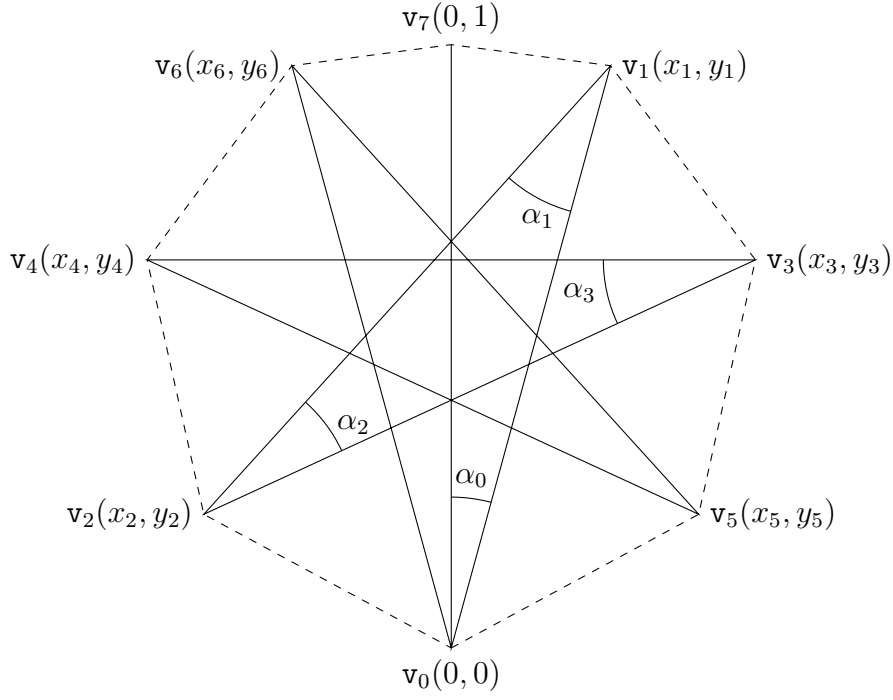

For all $n=2m$ with integer $m\ge 3$, consider a small $n$-gon $\geo{P}_n$ having the optimal diameter graph: an $(n-1)$-length cycle $\geo{v}_{0} - \geo{v}_1 - \ldots-\geo{v}_{k} - \ldots - \geo{v}_{\frac{n}{2}-1} - \geo{v}_{\frac{n}{2}} - \ldots - \geo{v}_{n-k-1} - \ldots - \geo{v}_{n-2}-\geo{v}_0$ plus the pendant edge $\geo{v}_{0} - \geo{v}_{n-1}$, as illustrated in Figure~\ref{figure:model:optimal}. We assume that $\geo{P}_n$ has the edge $\geo{v}_{0}-\geo{v}_{n-1}$ as axis of symmetry.

We use cartesian coordinates to describe the $n$-gon $\geo{P}_n$, assuming that a vertex $\geo{v}_k$, $k=0,1,\ldots,n-1$, is positioned at abscissa $x_k$ and ordinate $y_k$. Placing the vertex $\geo{v}_0$ at the origin, we set $x_0 = y_0 = 0$. We also assume that $\geo{P}_n$ is in the half-plane $y\ge 0$.

Let us place the vertex $\geo{v}_{n-1}$ at $(0,1)$ in the plane. Let $\alpha_0 = \angle \geo{v}_{n-1}\geo{v}_{0}\geo{v}_{1}$ and for all $k=1,2,\ldots, n/2-1$, $\alpha_k = \angle \geo{v}_{k-1} \geo{v}_{k} \geo{v}_{k+1}$. Since $\geo{P}_n$ is symmetric, we have
\begin{equation}\label{eq:condition}
\sum_{k=0}^{n/2-1}\alpha_k = \frac{\pi}{2},
\end{equation}
and
\begin{subequations}\label{eq:xy}
\begin{align}
x_{k} &= \sum_{i=0}^{k-1} (-1)^{i} \sin \left(\sum_{j=0}^{i}\alpha_j\right) && = - x_{n-k-1} &\forall k=1,2,\ldots, \frac{n}{2}-1,\\
y_{k} &= \sum_{i=0}^{k-1} (-1)^{i} \cos \left(\sum_{j=0}^{i}\alpha_j\right) && = y_{n-k-1} &\forall k=1,2,\ldots,\frac{n}{2}-1.
\end{align}
\end{subequations}
Since the edge $\geo{v}_{\frac{n}{2}-1} - \geo{v}_{\frac{n}{2}}$ is horizontal and $\|\geo{v}_{\frac{n}{2}-1} - \geo{v}_{\frac{n}{2}}\| = 1$, we also have
\begin{equation}\label{eq:x:m}
	x_{\frac{n}{2}-1} = (-1)^{\frac{n}{2}}/2 = -x_{\frac{n}{2}}.
\end{equation}

If $A_1$ denote the area of the triangle $\geo{v}_0 \geo{v}_1 \geo{v}_{n-1}$ and $A_k$ the area of the triangle $\geo{v}_0 \geo{v}_{k+1} \geo{v}_{k-1}$ for all $k = 2,3,\ldots,n/2-1$, then the area of $\geo{P}_n$ is $A = \sum_{k=1}^{n/2-1} 2A_k$. From~\eqref{eq:condition} and \eqref{eq:xy}, we have
\begin{subequations}\label{eq:A}
	\begin{align}
		2A_1 &= x_1 = \sin \alpha_0,\\
		2A_k &= x_{k+1}y_{k-1} - y_{k+1}x_{k-1} \nonumber\\
		&= \sin \alpha_k + 2(-1)^k \left(x_k \sin \left(\sum_{j=0}^{k-1} \alpha_{j} + \frac{\alpha_k}{2}\right)+ y_k \cos \left(\sum_{j=0}^{k-1} \alpha_{j} + \frac{\alpha_k}{2}\right) \right)\sin \frac{\alpha_k}{2}
	\end{align}
\end{subequations}
for all $k = 2,3,\ldots,n/2-1$. Then one can construct a large small $n$-gon by maximizing the area $A$ over $n/2$ variables $\alpha_0, \alpha_1, \ldots, \alpha_{\frac{n}{2}-1}$ subject to~\eqref{eq:condition} and~\eqref{eq:x:m}. Instead, we are going to use the same approach as Mossinghoff~\cite{mossinghoff2006b} to obtain a large small $n$-gon with fewer variables.

Now, suppose $\alpha_0 = \alpha$, $\alpha_1 = \beta + \gamma$, $\alpha_2 = \beta - \gamma$, and $\alpha_k = \beta$ for all $k = 3,4,\ldots, n/2-1$. Then \eqref{eq:condition} becomes
\begin{equation}\label{eq:condition:ab}
	\alpha + \left(\frac{n}{2} - 1\right)\beta = \frac{\pi}{2}.
\end{equation}
Coordinates $(x_k,y_k)$ in~\eqref{eq:xy} are given by
\begin{subequations}\label{eq:xyab}
	\begin{align}
		x_{1} &= \sin \alpha,\\
		y_{1} &= \cos \alpha,\\
		x_{2} &= \sin \alpha - \sin (\alpha+\beta+\gamma), \label{eq:xyab:2}\\
		y_{2} &= \cos \alpha - \cos (\alpha+\beta+\gamma),\\
	x_{k} &= x_{2} + \sum_{j=3}^k (-1)^{j-1} \sin (\alpha + (j-1)\beta) \nonumber\\
	&= x_{2} + \frac{\sin \left(\alpha + 3\frac{\beta}{2}\right) - (-1)^k\sin \left(\alpha + (2k-1)\frac{\beta}{2}\right)}{2\cos \frac{\beta}{2}} &\forall k=3,4,\ldots,\frac{n}{2}-1, \label{eq:xyab:k}\\
	y_{k} &= y_{2} + \sum_{j=3}^k (-1)^{j-1} \cos (\alpha + (j-1)\beta) \nonumber\\
	&= y_{2} + \frac{ \cos \left(\alpha + 3\frac{\beta}{2}\right) - (-1)^k\cos \left(\alpha + (2k-1)\frac{\beta}{2}\right)}{2\cos \frac{\beta}{2}} &\forall k=3,4,\ldots,\frac{n}{2}-1.
\end{align}
\end{subequations}
From \eqref{eq:x:m}, \eqref{eq:xyab:2}, \eqref{eq:xyab:k}, and \eqref{eq:condition:ab}, we deduce that
\begin{equation}\label{eq:condition:abc}
\sin (\alpha + \beta + \gamma) = \sin \alpha + \frac{\sin \left(\alpha + 3\frac{\beta}{2}\right)}{2\cos \frac{\beta}{2}}.
\end{equation}

The areas $A_k$ in~\eqref{eq:A} determined by $\alpha$, $\beta$, and $\gamma$ are
\[
\begin{aligned}
	2A_1 &= \sin \alpha,\\
	2A_2 &= \sin (2\beta) - \sin (\beta + \gamma),\\
	2A_k &= \sin \beta + 2(-1)^k \left(x_k \sin \left(\alpha + (2k-1)\frac{\beta}{2}\right)+ y_k \cos \left(\alpha + (2k-1)\frac{\beta}{2}\right) \right) \sin \frac{\beta}{2}\\
	&= \sin \beta - \tan \frac{\beta}{2} + 2(-1)^{k-1}  \left( 2\sin \frac{\beta + \gamma}{2} \sin \left((k-1)\beta - \frac{\gamma}{2}\right) - \frac{\cos ((k-2)\beta)}{2\cos \frac{\beta}{2}}\right)\sin \frac{\beta}{2}
\end{aligned}
\]
for all $k = 3,4,\ldots,n/2-1$. Using~\eqref{eq:condition:abc}, it follows that
\[
\sum_{k=3}^{n/2-1} 2A_k = \left(\frac{n}{2}-3\right)\left(\sin \beta - \tan \frac{\beta}{2}\right) + \left(\cos (\beta - \gamma) - \cos (2\beta) - \frac{1}{2}\right) \tan \frac{\beta}{2}.
\]
Thus,
\begin{equation}\label{eq:A:abc}
	\begin{aligned}
		A &= \sin \alpha + \sin (2\beta) - \sin (\beta + \gamma) \\
		&+ \left(\frac{n}{2}-3\right)\left(\sin \beta - \tan \frac{\beta}{2}\right) + \left(\cos (\beta - \gamma) - \cos (2\beta) - \frac{1}{2}\right) \tan \frac{\beta}{2}.
	\end{aligned}
\end{equation}
Note that, for $n=6$, we have $A = \sin \alpha + \sin (2\beta) - \sin (\beta + \gamma)$.

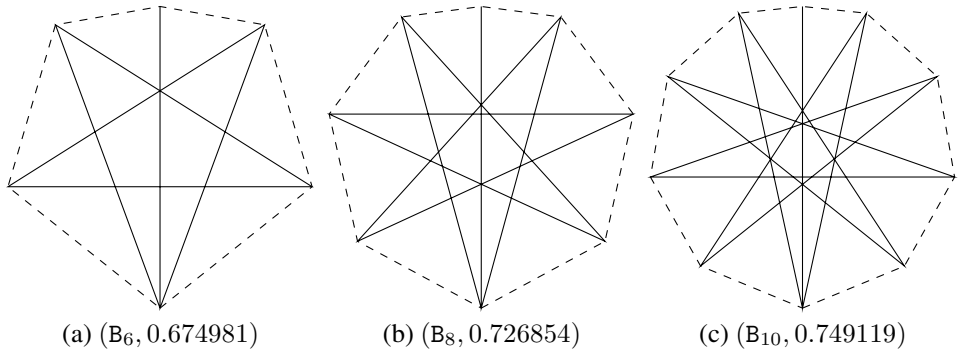
\begin{figure}[h]
	\centering
	\subfloat[$(\geo{B}_6,0.674981)$]{
		\begin{tikzpicture}[scale=4]
			\draw[dashed] (0,0) -- (0.5000,0.4024) -- (0.3438,0.9391) -- (0,1) -- (-0.3438,0.9391) -- (-0.5000,0.4024) -- cycle;
			\draw (0,1) -- (0,0) -- (0.3438,0.9391) -- (-0.5000,0.4024) -- (0.5000,0.4024) -- (-0.3438,0.9391) -- (0,0);
		\end{tikzpicture}
	}
	\subfloat[$(\geo{B}_{8},0.726854)$]{
		\begin{tikzpicture}[scale=4]
			\draw[dashed] (0,0) -- (0.4068,0.2215) -- (0.5000,0.6432) -- (0.2619,0.9651) -- (0,1) -- (-0.2619,0.9651) -- (-0.5000,0.6432) -- (-0.4068,0.2215) -- cycle;
			\draw (0,1) -- (0,0) -- (0.2619,0.9651) -- (-0.4068,0.2215) -- (0.5000,0.6432) -- (-0.5000,0.6432) -- (0.4068,0.2215) -- (-0.2619,0.9651) -- (0,0);
		\end{tikzpicture}
	}
	\subfloat[$(\geo{B}_{10},0.749119)$]{
		\begin{tikzpicture}[scale=4]
			\draw[dashed] (0,0) -- (0.3351,0.1395) -- (0.5000,0.4346) -- (0.4428,0.7678) -- (0.2103,0.9776) -- (0,1) -- (-0.2103,0.9776)  -- (-0.4428,0.7678) -- (-0.5000,0.4346) -- (-0.3351,0.1395) -- cycle;
			\draw (0,1) -- (0,0) -- (0.2103,0.9776) -- (-0.3351,0.1395) -- (0.4428,0.7678) -- (-0.5000,0.4346) -- (0.5000,0.4346) -- (-0.4428,0.7678) -- (0.3351,0.1395) -- (-0.2103,0.9776) -- (0,0);
		\end{tikzpicture}
	}
	\caption{Polygons $(\geo{B}_n,A(\geo{B}_n))$ defined in Theorem~\ref{thm:Bn}}
	\label{figure:Bn}
\end{figure}

With~\eqref{eq:condition:ab} and~\eqref{eq:condition:abc}, the area $A$ in~\eqref{eq:A:abc} can be considered as a one-variable function $f(\alpha)$. For instance, for $\alpha = \frac{\pi}{2n-2}$, we have $\beta = \frac{\pi}{n-1}$, $\gamma = 0$, and $f\left(\frac{\pi}{2n-2}\right) = A(\geo{R}_{n-1}^+)$. We may now search for a value of $\alpha \in \left[\frac{\pi}{2n-2}, \frac{\pi}{n}\right]$ that maximizes this function. An asymptotic analysis produces that, for large $n$, $f(\alpha)$ is maximized at $\hat{\alpha}(n)$ satisfying
\[
\hat{\alpha}(n) = \frac{a\pi}{n} + \frac{b\pi}{n^2} - \frac{c\pi}{n^3} + O\left(\frac{1}{n^4}\right),
\]
where $a = \frac{2\sqrt{114}-7}{22} = 0.652461\ldots$, $b = \frac{84a^2-272a+175}{4(22a+7)} = \frac{3521\sqrt{114}-34010}{9196} = 0.389733\ldots$, and
\[
\begin{aligned}
	c &= \frac{(7792a^4+16096a^3 + 2568a^2 -6248a +223)\pi^2}{768(22a+7)} - \frac{226a^2 + 84ab-22b^2-542a-136b + 303}{2(22a+7)}\\
	&= \frac{17328(663157+3161\pi^2) - (1088031703 - 3918085\pi^2)\sqrt{114}}{507398496} = 1.631188\ldots.
\end{aligned}
\]
Let $\geo{B}_n$ denote the $n$-gon obtained by setting $\alpha = \hat{\alpha}(n)$. We have
\[
\begin{aligned}
	\beta &= \hat{\beta}(n) = \frac{\pi}{n} + \frac{2(1-a)\pi}{n^2} + O\left(\frac{1}{n^3}\right),\\
	\gamma &= \hat{\gamma}(n) = \frac{(2a-1)\pi}{4n} + \frac{(a+b-1)\pi}{2n^2} + O\left(\frac{1}{n^3}\right),
\end{aligned}
\]
and the area of $\geo{B}_n$ is
\[
\begin{aligned}
A(\geo{B}_n) &= f(\hat{\alpha}(n))\\
&= \frac{\pi}{4} - \frac{5\pi^3}{48n^2} - \frac{(5545-456\sqrt{114})\pi^3}{5808n^3} - \left(\frac{7(13817-1281\sqrt{114})}{10648} - \frac{\pi^2}{480}\right) \frac{\pi^3}{n^4}\\
&- \left( \frac{23\pi^2(351468\sqrt{114}-2868731)}{618435840} + \frac{4013754104-375661161\sqrt{114}}{53410368}\right) \frac{\pi^3}{n^5} + O\left(\frac{1}{n^6}\right),
\end{aligned}
\]
which implies
\[
	\ub{A}_n - A(\geo{B}_n) = \frac{(5303-456\sqrt{114})\pi^3}{5808n^3} + \frac{(192107-17934\sqrt{114})\pi^3}{21296n^4} + O\left(\frac{1}{n^5}\right).
\]
By construction, $\geo{B}_n$ is small. We illustrate $\geo{B}_n$ for some $n$ in Figure~\ref{figure:Bn}.

Mossinghoff's small $n$-gons $\geo{M}_n$, $n=2m$ and $m\ge 3$, constructed in~\cite{mossinghoff2006b} for the maximal area problem were obtained as follows. He first supposed that $\alpha_0 = \alpha$, $\alpha_1 = \beta + \gamma$, $\alpha_2 = \beta - \gamma$, and $\alpha_k = \beta$ for all $k = 3,4,\ldots, n/2-3$. Then he set $\alpha = \frac{a\pi}{n} + \frac{t\pi}{n^2}$, $\beta = \frac{\pi}{n} + \frac{2(1-a)\pi}{n^2}$, and $\gamma = \frac{(2a-1)\pi}{4n} + \frac{(a+t-1)\pi}{2n^2}$, with
\[
\begin{aligned}
t &= \frac{4(7a^2-32a+25)}{44a+27} + (-1)^{\frac{n}{2}} \frac{15\pi(8a^3+12a^2-2a-3)}{32(44a+27)}\\
&= \frac{103104\sqrt{114}-998743}{200255} + (-1)^{\frac{n}{2}} \frac{15\pi(347\sqrt{114}-714)}{1762244}\\
&=
\begin{cases}
	0.429901\ldots &\text{if $n \equiv 2 \bmod 4$,}\\
	0.589862\ldots &\text{if $n \equiv 0 \bmod 4$.}
\end{cases}
\end{aligned}
\]
Note that we do not require $\alpha_{\frac{n}{2}-2} = \alpha_{\frac{n}{2}-1} = \beta$ in $\geo{M}_n$. The area of $\geo{M}_n$ is given by
\[
\begin{aligned}
	A(\geo{M}_n) &= \frac{\pi}{4} - \frac{5\pi^3}{48n^2} - \frac{(5545-456\sqrt{114})\pi^3}{5808n^3} - \left(\frac{7(13817-1281\sqrt{114})}{10648} - \frac{\pi^2}{480}\right) \frac{\pi^3}{n^4}\\
	&-\left(\frac{\pi^2(28622156724\sqrt{114}-177320884133)}{2251724893440} + \frac{182558364974-17072673147\sqrt{114}}{2326162080}\right.\\
	&+ (-1)^\frac{n}{2}\left. \frac{45\pi(1012477-919131\sqrt{114})}{852926096}\right) \frac{\pi^3}{n^5} + O\left(\frac{1}{n^6}\right),
\end{aligned}
\]
Therefore,
\[
A(\geo{B}_n) - A(\geo{M}_n) = \frac{3d\pi^3}{n^5} + O\left(\frac{1}{n^6}\right)
\]
with
\[
\begin{aligned}
	d &= \frac{25\pi^2(1747646-22523\sqrt{114})}{4691093528} + \frac{32717202988-3004706459\sqrt{114}}{29464719680}\\
	&+ (-1)^{\frac{n}{2}} \frac{15\pi(10124777-919131\sqrt{114})}{852926096}\\
	&=
	\begin{cases}
		0.0836582354\ldots &\text{if $n \equiv 2 \bmod 4$,}\\
		0.1180393778\ldots &\text{if $n \equiv 0 \bmod 4$.}
	\end{cases}
\end{aligned}
\]
We can also note that, for some parameter $u$,
\[
A(\geo{B}_n) - f\left(\frac{a\pi}{n} + \frac{u\pi}{n^2}\right) =
\begin{cases}
	\frac{(u-b)^2\pi^3\sqrt{114}}{8n^5} + O\left(\frac{1}{n^6}\right) &\text{if $u \not=b$,}\\
	\frac{c^2\pi^3\sqrt{114}}{8n^7} + O\left(\frac{1}{n^8}\right) &\text{if $u =b$.}
\end{cases}
\]
This completes the proof of Theorem~\ref{thm:Bn}.\qed

Table~\ref{table:area} shows the areas of $\geo{B}_n$, along with the optimal values $\hat{\alpha}(n)$, the upper bounds $\ub{A}_n$, the areas of $\geo{R}_n$, $\geo{R}_{n-1}^+$, and $\geo{M}_n$ for $n=2m$ and $3 \le m \le 12$. We also report the areas of the small $n$-gons $\geo{M}_n'$ obtained by setting $\alpha = \frac{a\pi}{n} + \frac{t\pi}{n^2}$ in~\eqref{eq:A:abc}, i.e., $A(\geo{M}_n') = f\left(\frac{a\pi}{n} + \frac{t\pi}{n^2}\right)$. Values in the table are rounded at the last printed digit. As suggested by Theorem~\ref{thm:Bn}, when $n$ is even, $\geo{B}_n$ provides a tighter lower bound on the maximal area $A_n^*$ compared to the best prior small $n$-gon~$\geo{M}_n$. For instance, we can note that $A(\geo{B}_{6}) = A_6^*$. We also remark that $A(\geo{M}_n) < A(\geo{M}_n')$ for all even $n\ge 8$.

\begin{table}[h]
	\footnotesize
	\centering
	\caption{Areas of $\geo{B}_n$}
	\label{table:area}
		\resizebox{\linewidth}{!}{
	\begin{tabular}{@{}rlllllll@{}}
		\toprule
		$n$ & $\hat{\alpha}(n)$ & $A(\geo{R}_n)$ & $A(\geo{R}_{n-1}^+)$ & $A(\geo{M}_n)$ & $A(\geo{M}_n')$ & $A(\geo{B}_n)$ & $\ub{A}_n$ \\
		\midrule
		6	&	0.3509301889	&	0.6495190528	&	0.6722882584	&	0.6731855653	&	0.6731855653	&	0.6749814429	&	0.6877007594	\\
		8	&	0.2649613582	&	0.7071067812	&	0.7253199909	&	0.7259763468	&	0.7264449921	&	0.7268542719	&	0.7318815691	\\
		10	&	0.2119285702	&	0.7347315654	&	0.7482573378	&	0.7490291363	&	0.7490910913	&	0.7491189262	&	0.7516135587	\\
		12	&	0.1762667716	&	0.7500000000	&	0.7601970055	&	0.7606471438	&	0.7606885130	&	0.7607153082	&	0.7621336536	\\
		14	&	0.1507443724	&	0.7592965435	&	0.7671877750	&	0.7675035228	&	0.7675178190	&	0.7675203660	&	0.7684036467	\\
		16	&	0.1316139556	&	0.7653668647	&	0.7716285345	&	0.7718386481	&	0.7718489998	&	0.7718535572	&	0.7724408116	\\
		18	&	0.1167583322	&	0.7695453225	&	0.7746235089	&	0.7747776809	&	0.7747819422	&	0.7747824059	&	0.7751926059	\\
		20	&	0.1048968391	&	0.7725424859	&	0.7767382147	&	0.7768497848	&	0.7768531741	&	0.7768543958	&	0.7771522071	\\
		22	&	0.0952114547	&	0.7747645313	&	0.7782865351	&	0.7783722564	&	0.7783738385	&	0.7783739622	&	0.7785970008	\\
		24	&	0.0871560675	&	0.7764571353	&	0.7794540033	&	0.7795196190	&	0.7795209668	&	0.7795213955	&	0.7796927566	\\
		\bottomrule
	\end{tabular}
		}
\end{table}

All polygons presented in this work and in~\cite{bingane2021a,bingane2021b,bingane2021c,bingane2021d,bingane2021e} were implemented as a MATLAB package: OPTIGON~\cite{optigon}, which is freely available at \url{https://github.com/cbingane/optigon}. In OPTIGON, we provide MATLAB functions that give the coordinates of their vertices. One can also find an algorithm developed in~\cite{bingane2021a} to find an estimate of the maximal area of a small $n$-gon when $n \ge 6$ is even.

\section{Conclusion}\label{sec:conclusion}
Tighter lower bounds on the maximal area of small $n$-gons were provided when $n$ is even. For each $n=2m$ with integer $m\ge 3$, we constructed a small $n$-gon $\geo{B}_n$ whose area is the maximum value of a one-variable function. For all even $n\ge 6$, the area of $\geo{B}_n$ is greater than that of the best prior small $n$-gon constructed by Mossinghoff. Furthermore, for $n=6$, $\geo{B}_6$ is the largest small $6$-gon.

\section*{Acknowledgements}
The author thanks Charles Audet, Professor at Polytechnique Montr\'{e}al, for helpful discussions on extremal small polygons and helpful comments on early drafts of this paper.

\bibliographystyle{ieeetr}
\bibliography{../../research}

\begin{thebibliography}{10}

\bibitem{reinhardt1922}
K.~Reinhardt, ``{Extremale polygone gegebenen durchmessers},'' {\em
  Jahresbericht der Deutschen Mathematiker-Vereinigung}, vol.~31, pp.~251--270,
  1922.

\bibitem{bieri1961}
H.~Bieri, ``{Ungel{\"o}ste probleme: Zweiter nachtrag zu nr. 12},'' {\em Elem.
  Math}, vol.~16, pp.~105--106, 1961.

\bibitem{graham1975}
R.~L. Graham, ``{The largest small hexagon},'' {\em Journal of Combinatorial
  Theory, Series A}, vol.~18, no.~2, pp.~165--170, 1975.

\bibitem{audet2002}
C.~Audet, P.~Hansen, F.~Messine, and J.~Xiong, ``{The largest small octagon},''
  {\em Journal of Combinatorial Theory, Series A}, vol.~98, no.~1, pp.~46--59,
  2002.

\bibitem{henrion2013}
D.~Henrion and F.~Messine, ``{Finding largest small polygons with
  GloptiPoly},'' {\em Journal of Global Optimization}, vol.~56, no.~3,
  pp.~1017--1028, 2013.

\bibitem{audet2017}
C.~Audet, ``{Maximal area of equilateral small polygons},'' {\em The American
  Mathematical Monthly}, vol.~124, no.~2, pp.~175--178, 2017.

\bibitem{foster2007}
J.~Foster and T.~Szabo, ``{Diameter graphs of polygons and the proof of a
  conjecture of Graham},'' {\em Journal of Combinatorial Theory, Series A},
  vol.~114, no.~8, pp.~1515--1525, 2007.

\bibitem{mossinghoff2006b}
M.~J. Mossinghoff, ``{Isodiametric problems for polygons},'' {\em Discrete \&
  Computational Geometry}, vol.~36, no.~2, pp.~363--379, 2006.

\bibitem{audet2009}
C.~Audet, P.~Hansen, and F.~Messine, ``{Ranking small regular polygons by area
  and by perimeter},'' {\em Journal of Applied and Industrial Mathematics},
  vol.~3, no.~1, pp.~21--27, 2009.

\bibitem{yuan2004}
B.~Yuan, ``{The Largest Small Hexagon},'' Master's thesis, National University
  of Singapore, 2004.

\bibitem{bingane2021a}
C.~Bingane, ``{Largest small polygons: A sequential convex optimization
  approach},'' Tech. Rep. G-2020-50, Les cahiers du GERAD, 2020.
\newblock \url{https://arxiv.org/abs/2009.07893}.

\bibitem{bingane2021b}
C.~Bingane, ``{Tight bounds on the maximal perimeter and the maximal width of
  convex small polygons},'' Tech. Rep. G-2020-53, Les cahiers du GERAD, 2020.
\newblock \url{https://arxiv.org/abs/2010.02490}.

\bibitem{bingane2021c}
C.~Bingane and C.~Audet, ``{Tight bounds on the maximal perimeter of convex
  equilateral small polygons},'' Tech. Rep. G-2021-31, Les cahiers du GERAD,
  2021.
\newblock \url{https://arxiv.org/abs/2105.10618}.

\bibitem{bingane2021d}
C.~Bingane, ``{Maximal perimeter and maximal width of a convex small
  polygon},'' Tech. Rep. G-2021-33, Les cahiers du GERAD, 2021.
\newblock \url{https://arxiv.org/abs/2106.11831}.

\bibitem{bingane2021e}
C.~Bingane and C.~Audet, ``{The equilateral small octagon of maximal width},''
  Tech. Rep. G-2021-45, Les cahiers du GERAD, 2021.
\newblock \url{https://arxiv.org/abs/2110.00036}.

\bibitem{optigon}
C.~Bingane, ``{OPTIGON: Extremal small polygons}.''
  \url{https://github.com/cbingane/optigon}, September 2020.

\end{thebibliography}
\end{document}